\documentclass[11pt]{article}
\usepackage{amssymb}
\usepackage{amsmath}
\usepackage{amsthm}
\usepackage{latexsym}
\usepackage{amsfonts}
\usepackage{graphicx}
\usepackage{graphics}

\newcommand{\dis}{\displaystyle}
\theoremstyle{plain}
\newtheorem{thm}{Theorem}[section]   % Αρίθμηση συνεχόμενη (όχι κατά Θεώρημα, Λήμμα κ.λπ.)

\newtheorem{lem}[thm]{Lemma}

\theoremstyle{definition}
\newtheorem{rem}[thm]{Remark}

\newtheorem*{Proof}{Proof}

\newcommand{\ra}{\;\rightarrow\;}

\newcommand{\Ga} {{\varGamma}}

\newcommand{\de}{\delta }
\newcommand{\OO} {{\varOmega}}
\newcommand{\De} {{\varDelta}}
\newcommand{\e}{\varepsilon }

\newcommand{\la}{\lambda }

\newcommand{\ti}{\tau }

\newcommand{\C}{\mathbb{C}}

\newcommand{\N}{\mathbb{N}}

\newcommand{\A}{\mathbb{A}}

\newcommand{\ssum}{\sum\limits}

\newcommand{\oO}{\overline{\varOmega}}

\newcommand{\ld}{\ldots}

\newcommand{\sm}{\smallsetminus}

\newcommand{\qb}{$\quad\blacksquare$}
\begin{document}
\title{\bf Generic Approximation of functions\\ by their Pad\'{e} approximants, II}
\author{G. Fournodavlos}
\date{}
\maketitle
%
%\begin{abstract}
\noindent
%
%{\bf Abstract}.\medskip
\begin{abstract}

In \cite{5} we proved that generically functions defined in any open set can be approximated by a sequense of their Pad\'{e} approximants, in the sense of uniform convergence on compacta. In this paper we examine a more particular space, $A^{\infty}(\Omega)$ , and prove that we can obtain similar approximation results with functions smooth on the boundary.
\end{abstract}
{\em Subject Classification MSC2010}\,: primary 41A21, 30K05
secondary 30B10, 30E10, 30K99, 41A10, 41A20. \vspace*{0.2cm} \\
{\em Key words}\,: Pad\'{e} approximant, Taylor series, Baire's
theorem, Runge's theorem, generic property.
\section{Introduction}\label{sec1} % 1
\noindent

On a disc every holomorphic function $f$ can be approximated by the partial sums of its Taylor expansion. It is also true that generically in a simply connected domain every holomorphic function is the uniform on compacta limit of a subsequense of the partial sums of its Taylor expansion. The partial sums are polynomials and thus by the maximum principle we are led to uniform approximation on compact sets with connected complement. If we replace the partial sums by the Pad\'{e} approximants $[p/q]_f$, which are rational functions with poles, then we obtain approximation on compact sets with arbitrary connectivity (\cite{4}, \cite{5}).

In \cite{2} it was proved that generically every entire function can be approximated uniformly on compacta by a subsequense $[p_n/q_n]_f$ of its Pad\'{e} approximants, provided $p_n\ra+\infty$ and $p_n-q_n\ra+\infty$. In \cite{5} we weakened the previous assumption to $p_n\ra+\infty$ only and we extended the result to any simply connected domain. We also obtained the same approximation on any open subset of $\C$ (of arbitrary connectivity) under the assumption $p_n\ra+\infty$ and $q_n\ra+\infty$.

In the present paper we obtain similar results using smooth functions; that is, holomorphic functions on $\OO$ such that every derivative $f^{(l)}$ extends continuously on $\oO$ $(f\in A^{\infty}(\OO))$. In the case of a domain $\OO\subseteq\C$ such that $(\C\cup\{\infty\}\sm\oO$ is connected we obtain the result provided $p_n\ra+\infty$; this result is generic in a subset of $A^{\infty}(\OO)$, which is the closure of the set of polynomials, under the natural topology. We do not know in general if this subspace is the whole $A^{\infty}(\OO)$. If there exists a finite constant $M<\infty$ so that all points $A,B\in \oO$ can be joined in $\oO$ by a curve $\Ga$ with length $|\Ga|\le M$ , then the polynomials are dense in $A^{\infty}(\OO)$ (\cite{11}). In the case $q=0$ the generic result in the closure of polynomials in $A^{\infty}(\OO)$ is known (\cite{9}).

Finally in the general case of an open subset $\OO\subseteq\C$ we obtain a similar generic result in the closure in $A^{\infty}(\OO)$ of holomorphic functions in some varying neighborhood of $\OO$, provided $p_n\ra+\infty$ and $q_n\ra+\infty$. Our method of proof is based on Baire's Category theorem (\cite{7}, \cite{8}) and extends the methods of \cite{4} and \cite{13}.
\section{Preliminaries}\label{sec2}
\noindent

Let $\OO\subseteq\C$ be an open set and let us consider the set $A^{\infty}(\OO)=\{ f\in H(\OO)$: $f^{(l)}$ extends continuously on $\oO$, $l=0,1,\ld$\}, where $f^{(l)}, l=0,1,\ld$ denote the derivatives of the holomorphic function $f$.

We define the following metric $\rho$ on $A^{\infty}(\OO)$:
\[
\rho(f,g)=\sum^\infty_{l=0}\sum^\infty_{n=1}\frac{1}{2^{l+n}}\min\{\|f-g\|_{l,n},1\},
\]
where $\|f-g\|_{l,n}=\dis\sup_{z\in \oO\cap\overline{\De(0,n)}}|f^{(l)}-g^{(l)}|$, $l=0,1,\ld$. It is easy to see that a sequense in $A^{\infty}(\OO)$, $(f_m)_{m\in\N}$ converges $f_m\overset{\rho}{\longrightarrow}f\in A^{\infty}(\OO)$, if and only if ${f_m}^{(l)}\ra f^{(l)}$
uniformly on each compact subset of $\oO$, for every $l$. The space $(A^{\infty}(\OO),\rho)$
is compete.

Let $f$ be a function holomorphic in a neighborhood of 0 and let
$f(z)=\ssum^\infty_{v=0}a_vz^v$ its Taylor series. A Pad\'{e}
approximant $[p/q]_f$ of $f$, $p,q\in\{0,1,2,\ld\}$, is a rational
function of the form
\[
\frac{\ssum^p_{v=0}n_vz^v}{\ssum^q_{v=0}d_vz^v}, \ \ d_0=1.
\]
such that its Taylor series $\ssum^\infty_{v=0}b_vz^v$ coincides
with $\ssum^\infty_{v=0}a_vz^v$ up to the first $p+q+1$ terms;
that is $b_v=a_v$ for $v=0,\ld,p+q$ (\cite{1}).

We notice that in case of $q=0$ there exists always a unique
Pad\'{e} approximant of $f$ and $[p/q]_f(z)= S_p(z)$, where $S_p(z)= \dis\sum^p_{v=0}a_vz^v$.
For $q\ge 1$ it is true that there exists a unique Pad\'{e} approximant of $f$,
if and only if the following
determinant is not zero:
\[
\det\left|\begin{array}{cccc}
  a_{p-q+1} & a_{p-q+2} & \cdots & a_p \\
  a_{p-q+2} & a_{p-q+3} & \cdots & a_{p+1} \\
  \vdots & \vdots &  & \vdots \\
  a_p & a_{p+1} & \cdots & a_{p+q-1} \\
\end{array}\right|\neq0, \ \ a_i=0, \ \ \text{when} \ \ i<0.
\eqno(\ast)
\]
Then we write $f\in D_{p,q}$.

If $f\in D_{p,q}$, then $[p/q]_f$ $(q\ge1)$ is given by the Jacobi
explicit formula:
\[
[p/q]_f=\frac{\det\left|\begin{array}{cccc}
  z^qS_{p-q}(z) & z^{q-1}S_{p-q+1} & \cdots & S_p(z) \\
  a_{p-q+1} & a_{p-q+2} & \cdots & a_{p+1} \\
  \vdots & \vdots &  & \vdots \\
  a_p & a_{p+1} & \cdots & a_{p+q} \\
\end{array}\right|}{\det\left|\begin{array}{cccc}
  z^q & z^{q-1} & \cdots & 1 \\
  a_{p-q+1} & a_{p-q+2} & \cdots & a_{p+1} \\
  \vdots & \vdots &  & \vdots \\
  a_p & a_{p+1} & \cdots & a_{p+q} \\
\end{array}\right|},
\]
with $S_k(z)=\left\{\begin{array}{cc}
  \ssum^k_{v=0}a_vz^v, & k\ge0 \\
  0, & k<0. \\
\end{array}\right.$
\begin{rem}\label{rem2.1}
If all of the coefficients $\frac{f^{(v)}(0)}{v!}=a_v$,
$v=0,1,\ld,p+q$, involved in the determinant $(\ast)$ depend
linearly on $d\in\C$, $a_v=c_v\cdot d+\ti_v$, such that $c_v=0$,
when $v<p$ and $c_p\neq0$, then the determinant is a
polynomial in $d$ of degree $q$ and hence only for finite values of
$d$ the determinant is zero.
\end{rem}
If $L$ is any set we write $h\in H(L)$ if $h$ is holomorphic in
some open set containing $L$. We also denote $\|h\|_L=\dis\sup_{z\in L}|h(z)|$, for every function $h:L\rightarrow\C$ on the set $L$.
\begin{lem}\label{lem2.2}
Let $r>0$, $p,q,s\in\N$ and $K\subseteq\C$ a compact set. If $f\in\ H(\overline{\De(0,r)})$, $f\in D_{p,q}$ such that its Pad\'{e} approximant $[p/q]_f$ has no poles in $K$, then for every $\e>0$ there exists $\de>0$ such that for every $g\in\ H(\overline{\De(0,r)})$ with $\|g-f\|_{\overline{\De(0,r)}}<\de$ it holds $g\in\ D_{p,q}$ and $\|[p/q]^{(l)}_g-[p/q]^{(l)}_f\|_K<\e$, $\forall l\in\{0,1,\ld,s\}$.
\end{lem}
\begin{Proof}
Let $\e>0$. Observe that the determinant $(\ast)$ and the coefficients of the numerator and the denominator of $[p/q]^{(l)}_f$, $l=0,1,\ld,s$ depend polynomially on $\frac{f^{(v)}(0)}{v!}$, $v=0,1,\ld,p+q$. This implies that there exists $\widetilde{\de}$ such that for every $g\in H(\overline{\De(0,r)})$ with $\big|\frac{g^{(v)}(0)}{v!}-\frac{f^{(v)}(0)}{v!}\big|<\widetilde{\de}$, $v=0,1,\ld,p+q$ it holds $g\in D_{p,q}$ and $\|[p/q]^{(l)}_g-[p/q]^{(l)}_f\|_K<\e$, $l=0,1,\ld,s$.

If $0<\de<\min\{r^v\cdot\widetilde{\de}\;|\; v=0,1,\ld,p+q\}$ and $\|g-f\|_{\overline{\De(0,r)}}<\de$, then by Cauchy's estimates we obtain:
\[
\bigg|\frac{g^{(v)}(0)}{v!}-\frac{f^{(v)}(0)}{v!}\bigg|=
\bigg|\frac{(g-f)^{(v)}(0)}{v!}\bigg|\le\frac{\|g-f\|_{\overline{\De(0,r)}}}{r^v}<\frac{\de}{r^v}<\widetilde{\de}.
\qquad \text{\qb}
\]
\end{Proof}
\begin{rem}\label{rem2.3}
It follows from Lemma \ref{lem2.2} that $D_{p,q}\cap\ A^{\infty}(\OO)$ is open ($0\in \OO$).
\end{rem}
\section{A special case}\label{sec3}
\noindent

Let $\OO\subseteq\C$ be an open set containing 0, such that $(\C\cup\infty)\sm\oO$ is connected.
Also, let $F\subseteq\N\times\N$ which contains a sequence
$(\widetilde{p}_m,\widetilde{q}_m)_{m\in\N}$, such that
$\widetilde{p}_m\ra+\infty$. We define
\begin{itemize}
\item $B_F=\{f\in A^{\infty}(\OO)$: there exists $(p_m,q_m)_{m\in\N}$ in
$F$ such that $f\in D_{p_m,q_m}$, for all $m\in\N$ and for every
$K\subseteq\oO$ compact $[p_m/q_m]^{(l)}_f\ra f^{(l)}$ uniformly on $K$, for each $l=0,1,\ld \}$.
\item $E(n,s,(p,q))=\{f\in A^{\infty}(\OO):f\in D_{p,q}$ and
$\|[p/q]_f-f\|_{l,n}<1/s, l=0,1,\ld,s\}$, $n,s\in\N, (p,q)\in F$.
\end{itemize}
\begin{lem}\label{lem3.1}
$B_F=\bigcap\limits^\infty_{n,s=1}\bigcup\limits_{(p,q)\in
F}E(n,s,(p,q)$.
\end{lem}
\begin{Proof}
It is standard and is omitted. [A similar proof can be found in
\cite{12}].
\qb
\end{Proof}
\begin{lem}\label{lem3.2}
$E(n,s,(p,q))$ is open.
\end{lem}
\begin{Proof}
$D_{p,q}\cap\A^{\infty}(\OO)$ is open (Remark \ref{rem2.3}) and similarly to the proof of the Lemma
\ref{lem2.2}, we can prove that the map $f\mapsto\|[p/q]_f-f\|_{l,n}$ is continuous, for any $l$.
\qb
\end{Proof}

We will now focus our attention on a more accessible space, $H(\oO)$, which is a subspace of $A^{\infty}(\OO)$
and is considered with its relative topology.
\begin{lem}\label{lem3.3}
The polynomials are dense in $H(\oO)$.
\end{lem}
\begin{Proof}
Let $f\in H(\oO)$ and $\e>0$. It suffices to show that for $N=N(\e)\in\N$ and $L=L(\e)\in\N$
there exists a polynomial $P$ such that $\|P-f\|_{l,N}<\e$, $\forall l\le L$. Observe that
$(\C\cup\infty)\sm(\oO\cap\overline{\De(0,N)})$ is connected.
\begin{itemize}
\item $f\in H(\oO)$, thus, there exists $U\subseteq\C$ open such that $f\in H(U)$ and $\oO\subseteq U$.
\item It is true that we can find $V\subseteq\C$ open, such that $\oO\cap\overline{\De(0,N)}\subseteq V\subseteq U$
and $(\C\cup\infty)\sm V$ connected (in other words $V$ is simply connected) (\cite{3}, \cite{6}).
\end{itemize}

By Runge's theorem there exists a sequense of polynomials $(P_i)_{i\in\N}$,
such that $P^{(l)}_i\ra f^{(l)}$ uniformly on each compact subset of $V$
for every $l$, which completes the proof.
\qb
\end{Proof}
\begin{thm}\label{thm3.4}
$B_F\cap cl_{A^{\infty}(\OO)}H(\oO)$ is $G_\de$ and dense in $cl_{A^{\infty}(\OO)}H(\oO)$. (Hence $B_F\neq\emptyset$).
\end{thm}
\begin{Proof}
Lemma \ref{lem3.2} implies that $\bigcup\limits_{(p,q)\in
F}E(n,s(p,q))\cap cl_{A^{\infty}(\OO)}H(\oO)$ is open in $cl_{A^{\infty}(\OO)}H(\oO)$. By Lemma \ref{lem3.1} $B_F\cap cl_{A^{\infty}(\OO)}H(\oO)$ is $G_\de$ in $cl_{A^{\infty}(\OO)}H(\oO)$. We
claim that $\bigcup\limits_{(p,q)\in F}E(n,s,(p,q))\cap cl_{A^{\infty}(\OO)}H(\oO)$ is dense in $cl_{A^{\infty}(\OO)}H(\oO)$. If
this is true, then Baire's Category theorem completes the proof.
By Lemma \ref{lem3.3} it suffices to prove that for every polynomial $P$ and
$\e>0$ there exists $f\in\bigcup\limits_{(p,q)\in F}E(n,s,(p,q))\cap H(\oO)$
such that $\|P-f\|_{l,N}<\e$, for every $l\le L=L(\e)\in\N$, where $N=N(\e)\in\N$.
\begin{itemize}
\item Let $P$ be a polynomial and $\e>0$. There exists $(p,q)\in F$ such that
$p>\text{deg}P$.
\end{itemize}

If $q=0$, define $f(z)=P(z)+dz^p$, $d\in\C\sm\{0\}$. It is
immediate that $f\in D_{p,q}$ and $[p/q]_f=f$. It follows $f\in
E(n,s,(p,q))\cap H(\oO)$. In addition,
$\|f-P\|_{l,N}=|d|\cdot\|z^p\|_{l,N}<\e$, $\forall l\le L$, when
$0<|d|<\e/\dis\max_{0\le l\le L}\|z^p\|_{l,N}$.

If $q\ge1$, we define
$\widetilde{f}_j(z)=\frac{P(z)+d_jz^p}{1-(c_jz)^q}$, $c_j,d_j\in\C\sm\{0\}$,
where $d_j$ and $c_j$ will be determined later on, $j\in\N$.
\begin{itemize}
\item Let $\la>\max\{n,N\}$. We have $\dis\inf_{z\in\overline{\De(0,\la)}}|1-(c_jz)^q|\ge1-|c_j|^q\cdot\|z\|^q_{\overline{\De(0,\la)}}>\frac{1}{2}$,
when $0<|c_j|<\frac{1}{2^{\frac{1}{q}}\|z\|_{\overline{\De(0,\la)}}}$, $j\in\N$.
\item We have
$\|\widetilde{f}_j(z)-P(z)\|_{\overline{\De(0,\la)}}=\big\|\frac{P(z)+d_jz^p}{1-(c_jz)^q}\big\|_{\overline{\De(0,\la)}}\le
2(\|P(z)\|_{\overline{\De(0,\la)}}\cdot|c_j|^q\cdot\|z\|^q_{\overline{\De(0,\la)}}+|d_j|\cdot\|z\|^q_{\overline{\De(0,\la)}})$. Thus, there exists $\de_j>0$, $j\in\N$, such that: $\|\widetilde{f}_j-P\|_{\overline{\De(0,\la)}}<1/j$, when $0<|c_j|<\de_j<
\frac{1}{2^{\frac{1}{q}}\|z\|_{\overline{\De(0,\la)}}}$ and $0<|d_j|<\de_j$. Hence, $\widetilde{f}_j\ra P$ uniformly on ${\overline{\De(0,\la)}}$ and so $\widetilde{f}^{(l)}_j\ra P^{(l)}$ uniformly on $\oO\cap\overline{\De(0,N)}$ (which is contained in $\De(0,\la)$), for every $l$. Therefore, there exists $j_0\in\N$ such that $\|\widetilde{f}_{j_0}-P\|_{l,N}<\e/2$, $l=0,1,\ld,L$.
\item We fix $c_{j_0}$ satisfying the above. Around 0,
$\widetilde{f}_{j_0}(z)=P(z)+d_{j_0}z^p+P(z)\cdot(c_{j_0}z)^q+d_{j_0}z^p\cdot(c_{j_0}z)^q+\cdots$.
According to Remark \ref{rem2.1} we can choose $0<|d_{j_0}|<\de_{j_0}$, such that $\widetilde{f}_{j_0}\in D_{p,q}$. By the uniqueness of the Pad\'{e} approximant of $\widetilde{f}_{j_0}$ we obtain $[p/q]_{\widetilde{f}_{j_0}}=\widetilde{f}_{j_0}$.
\item Let $r>0$: $\overline{\De(0,r)}\subseteq\OO\cap\De(0,\la)$.
By Lemma \ref{lem2.2} there exists $\de>0$ such that for every $f\in H(\overline{\De(0,r)})$ with $\|f-\widetilde{f}_{j_0}\|_{\overline{\De(0,r)}}<\de$
it holds $\|[p/q]_f-[p/q]_{\widetilde{f}_{j_0}}\|_{l,n}<1/2s$, $\forall l\le s$.
Also, we demand $0<\de<\min\{1/2s,\e/2\}$.
\item The Taylor series of $\widetilde{f}_{j_0}$ around 0 has radius of convergence greater than $\frac{1}{|c_{j_0}|}>\la$.
It follows that its partial sums
$\big(\sum^k_{v=0}\frac{\widetilde{f}^{(v)}_{j_0}(0)}{v!}z^v\big)^{(l)}\ra\widetilde{f}^{(l)}_{j_0}$
uniformly on $\oO\cap\overline{\De(0,\la)}$, for every $l$. Hence, there exists a partial sum
$f(z)=\sum^{k_0}_{v=0}\frac{\widetilde{f}^{(v)}_{j_0}(0)}{v!}z^v$ such that $\|f-\widetilde{f}_{j_0}\|_{l,\la}<\de$,
$\forall l\le\max\{s,L\}$.
\item $f$ satisfies:
$\|[p/q]_f-f\|_{l,n}\le\|[p/q]_f-[p/q]_{\widetilde{f}_{j_0}}\|_{l,n}
+\|\widetilde{f}_{j_0}-f\|_{l,n}<1/2s+\de<1/s$, $\forall l\le s$. It follows that
$f\in E(n,s,(p,q))\cap H(\oO)$. Also, it holds
$\|f-P\|_{l,N}\le\|f-\widetilde{f}_{j_0}\|_{l,N}+\|\widetilde{f}_{j_0}-P\|_{l,N}<
\e/2+\de<\e$, $\forall l\le L$.
\end{itemize}
This completes the proof.  \qb
\end{Proof}
\section{The general case}\label{sec4}
\noindent

Let $\OO\subseteq\C$ be an open set containing 0. Also, let
$F\subseteq\N\times\N$ which contains a sequence
$(\widetilde{p}_m,\widetilde{q}_m)_{m\in\N}$ such that
$\widetilde{p}_m\ra+\infty$ and $\widetilde{q}_m\ra+\infty$. We
define $B_F$ and $E(n,s,(p,q))$ similarly as in Section
\ref{sec3}.

The analogue of Lemmas \ref{lem3.1}, \ref{lem3.2} hold in this
case also. Like before we concentrate on $H(\oO)$ and its closure in $A^{\infty}(\OO)$.
\begin{lem}\label{lem4.1}
The rational functions with poles off $\oO$ are dense in $H(\oO)$.
\end{lem}
\begin{Proof}
Let $f\in H(\oO)$. There exists $U\subseteq\C$ open (depending on $f$) such that
$\oO\subseteq U$ and $f\in H(U)$. By Runge's theorem there exists a sequense $(R_i)_{i\in\N}$
of rational functions with poles in $(\C\cup\infty)\sm U$, hence $R_i\in H(U)\subseteq H(\oO)$, $\forall i\in\N$,
such that $R_i\ra f$ uniformly on each compact set of $U$. Similarly to the Lemma \ref{lem3.3},
for a given $\e$ there exists $i_0$ such that $\|R_i-f\|_{l,N}<\e/2$, $\forall l\le L$,
where $L=L(\e)\in\N$, $N=N(\e)\in\N$ are chosen so that $\rho(f,R_{i_0})<\e$. \qb
\end{Proof}
\begin{thm}\label{thm4.2}
$B_F\cap cl_{A^{\infty}(\OO)}H(\oO)$ is $G_\de$ and dense in $cl_{A^{\infty}(\OO)}H(\oO)$. (Hence $B_F\neq\emptyset$).
\end{thm}
\begin{Proof}
Since $\bigcup\limits_{(p,q)\in F}E(n,s,(p,q))$ is open, it follows
that $B_F\cap cl_{A^{\infty}(\OO)}H(\oO)$ is $G_\de$ in the subspace. By Baire's Category theorem the proof would be complete if the set $\bigcup\limits_{(p,q)\in F}E(n,s,(p,q))\cap H(\oO)$ was dense in $H(\oO)$, $n,s\in\N$.
By Lemma \ref{lem4.1} it suffices to show that for every rational function $R$ with poles off $\oO$ (or $R\in H(\oO)$)
and every $\e>0$ there exists $f\in\bigcup\limits_{(p,q)\in F}E(n,s,(p,q))\cap H(\oO)$ such that
$\|f-R\|_{l,N}<\e$, $\forall l\le L=L(\e)\in\N$, $N=N(\e)\in\N$.
\begin{itemize}
\item Let $R(z)=\frac{A(z)}{B(z)}$ be a rational function with
poles only in $(\C\cup\{\infty\})\sm\oO$, where $A$, $B$ are
polynomials and let $\e>0$. There exists $(p,q)\in F$ such that $p>\text{deg}A$ and
$q>\text{deg}B$. We define
$\widetilde{f}_j(z)=\frac{A(z)+d_jz^p}{B(z)-(c_jz)^q}$,
$c_j,d_j\in\C\sm\{0\}$, where $c_j$ and $d_j$ will be determined later on, $j\in\N$.
\item Since $R$ has no poles in $\oO$, there exists $U\subseteq\C$ open such that $R\in H(U)$ and $\oO\subseteq U$.
Also, there exists $K\subseteq U$ compact such that $K^0\supseteq\oO\cap\overline{\De(0,\la)}$, where $\la=\max\{n,N\}$,
and every component of $(\C\cup\infty)\sm K$ contains at least one component of $(\C\cup\infty)\sm U$ (\cite{14}).
\item We have $B(0)\neq 0$ and $\dis\inf_{z\in K}|B(z)|>0$. Furthermore,
$\dis\inf_{z\in K}|B(z)-(c_jz)^q|\ge\dis\inf_{z\in K}|B(z)|-|c_j|^q\cdot\|z\|^q_K>0$,
when $0<|c_j|<\big(\frac{\inf_K|B(z)|}{\|z\|^q_K}\big)^{1/q}$, $\forall j\in\N$. Thus,
\[ \|\widetilde{f}_j(z)-R(z)\|_K= \bigg\|\frac{A(z)(c_jz)^q+B(z)d_jz^p}{B(z)(B(z)-(c_jz)^q)}\bigg\|_K\le
\frac{\|A(z)\|_K\cdot|c_j|^q\cdot\|z\|^q_K+\|B(z)\|_K\cdot|d_j|\cdot\|z\|^p_K}{\dis\inf_{z\in K}|B(z)|\cdot\dis\inf_{z\in K}|B(z)-(c_jz)^q|}. \]
\item There exists $\de_j>0$, $j\in\N$, such that
$\|\widetilde{f}_j-R\|_K<1/j$, whenever $|c_j|<\de_j<\big(\frac{\inf_K|B(z)|}{\|z\|^q_K}\big)^{1/q}$ and
$|d_j|<\de_j$, $\forall j\in\N$. Hence, $\widetilde{f}_j\ra R$ uniformly on $K$ and
$\widetilde{f}^{(l)}_j\ra R^{(l)}$ on each compact subset of $K^0$, for every $l$.
This implies that there exists $j_0\in\N$, such that $\|\widetilde{f}_{j_0}-R\|_{l,N}<\e/2$, $\forall l\le L$.
\item We fix $c_{j_0}$ satisfying the above. Around 0 we have:
$\widetilde{f}_{j_0}(z)=B^{-1}(0)A(z)+B^{-1}(0)d_{j_0}z^p-B^{-1}(0)A(z)\cdot(\widetilde{B}(z)-1)-
B^{-1}(0)d_{j_0}z^p\cdot(\widetilde{B}(z)-1)+ \cdots$,
where $\widetilde{B}(z)=B^{-1}(0)B(z)-B^{-1}(0)(c_{j_0}z)^q$.
By Remark \ref{rem2.1} we can choose
$0<|d_{j_0}|<\de_{j_0}$ such that $\widetilde{f}_{j_0}\in D_{p,q}$.
Thus, there exists a unique Pad\'{e} approximant of $\widetilde{f}_{j_0}$ and
$\widetilde{f}_{j_0}=\frac{B^{-1}(0)A(z)+B^{-1}(0)d_{j_0}z^p}{B^{-1}(0)B(z)-B^{-1}(0)(c_{j_0}z)^q}$
satisfies $[p/q]_{\widetilde{f}_{j_0}}=\widetilde{f}_{j_0}$.
\item There exists $r>0$: $\overline{\De(0,r)}\subseteq\OO\cap\De(0,\la)\subseteq\ K^0$.
Lemma \ref{lem2.2} provides $0<\de<\min\{1/2s,\e/2\}$ such
that for every $f\in H(\overline{\De(0,r)})$ with
$\|f-\widetilde{f}_{j_0}\|_{\overline{\De(0,r)}}<\de$, it follows $f\in D_{p,q}$ and
$\|[p/q]_f-[p/q]_{\widetilde{f}_{j_0}}\|_{l,n}<1/2s$, $\forall l\le s$.
\item By Runge's theorem there exists a sequense of rational functions, $(R_i)_{i\in\N}$,
with poles off $K$ and more particularly (see previous property of $K$) off $U$, such that
$R_i\ra\widetilde{f}_{j_0}$ uniformly on $K$. This implies that $R^{(l)}_i\ra\widetilde{f}^{(l)}_{j_0}$
uniformly on each compact subset of $K^0$. Hence, there exists $f=R_{i_0}\in H(U)\subseteq H(\oO)$
such that $\|f-\widetilde{f}_{j_0}\|_{l,\la}<\de$, $\forall l\le\max\{s,L\}$,
because $\oO\cap\overline{\De(0,\la)}\subseteq K^0$.
\item It follows that
$\|[p/q]_f-f\|_{l,n}\le\|[p/q]_f-[p/q]_{\widetilde{f}_{j_0}}\|_{l,n}+\|\widetilde{f}_{j_0}-f\|_{l,n}<1/2s+\de<1/s$,
$\forall l\le s$.
Thus,
$f\in E(n,s,(p,q))\cap H(\oO)$. Moreover, it holds
$\|f-R\|_{l,N}\le\|f-\widetilde{f}_{j_0}\|_{l,N}+\|\widetilde{f}_{j_0}-R\|_{l,N}<\de+\e/2<\e$,
$\forall l\le L$.
\end{itemize}
This completes the proof. \qb
\end{Proof}
\bigskip
Department of Mathematics, \\
University of Athens \\
Panepistemiopolis, 157 84 Athens, Greece \\
e-mail: gregdavlos@hotmail.com


\begin{thebibliography}{14}
\bibitem{1} G. A. Baker, Jr. and P. R. Graves-Morris: Pad\'{e}
Approximants, Vol. 1 and 2, (Encyclopedia of Mathematics and its
Applications), Cambridge Un. Press, 2010.
%
\bibitem{2} P. B. Borwein: The usual behaviour of rational
approximants, Canad. Math. Bull. Vol. 26, 1983, p. 317-323.
%
\bibitem{3} G. Costakis: Some remarks on universal functions
and Taylor series, Math. Proc. Cambr. Philos. Soc., Vol. 128, 2000, p. 157-175.
%
\bibitem{4} N. J. Daras and V. Nestoridis: Universal Pad\'{e}
approximation, arXiv: 1102.4782v1, 2011.
%
\bibitem{5} G. Fournodavlos: Generic Approximation of functions
by their Pad\'{e} approximants I, arXiv: 1103.3841, 2011.
%
\bibitem{6} K.-G. Grosse-Erdmann: Holomorphe Monster und
universelle Funktionen, Mitt. Math. Sem. Giessen, Vol. 176, 1987, p. 1-84.
%
\bibitem{7} K.-G. Grosse-Erdmann: Universal families and
hypercyclic operators, Bull. Amer. math. Soc., Vol. 36, 1999, p.
345-381.
%
\bibitem{8} J.-P. Kahane: Baire's Category theorem and
trigonometric series, J. Anal. Math., Vol. 80, 2000, p. 143-182.
%
\bibitem{9} Ch. Kariofillis, Ch. Konstandaki and V. Nestoridis:
Smooth universal Taylor series, Monatsh. Math., Vol. 147, Issue 3, 2006, p. 249-257.
%
\bibitem{10} A. Melas and V. Nestoridis: Universality of Taylor
series as a generic property of holomorphic functions, Adv. Math.
Vol. 157, 2001, p. 138-176.
%
\bibitem{11} A. Melas and V. Nestoridis: On various types
of universal Taylor series, Complex Variables Theory, Vol. 44, Issue
3, 2001, p. 245-258.
%
\bibitem{12} V. Nestoridis: Universal Taylor series, Annales de l'
Institute Fourier, Vol. 46, 1996, p. 1293-1306.
%
\bibitem{13} V. Nestoridis: A strong notion of universal
Taylor series, J. London Math. Soc., Vol. 68, Issue 2, 2003, p. 712-724.
%
\bibitem{14} W. Rudin, Real and Complex Analysis, McGraw-Hill,
1986.
\end{thebibliography}
\end{document}